\documentclass[12pt]{amsart}
\usepackage{amscd}
\usepackage{amssymb}
\usepackage[all]{xy}

\numberwithin{equation}{section}

\swapnumbers

\theoremstyle{plain} 
\newtheorem{thm}[equation]{Theorem}
\newtheorem{cor}[equation]{Corollary}

\newtheorem{prop}[equation]{Proposition}
\newtheorem{conj}[equation]{Conjecture}

\theoremstyle{definition}
\newtheorem{defn}[equation]{Definition}

\theoremstyle{remark}
\newtheorem{rem}[equation]{Remark}
\newtheorem{rems}[equation]{Remarks}
\newtheorem{ex}[equation]{Example}
\newtheorem{problem}[equation]{Problem}

\newtheorem{exs}[equation]{Examples}

\newtheorem{exercise}[equation]{Exercise}

\newtheorem*{VariableNoNum}{{\VariableText}}
\newtheorem{Variable}[equation]{{\VariableText}}

\theoremstyle{definition}
\newtheorem*{VariableNoNumBold}{{\VariableText}}
\newtheorem{VariableBold}[equation]{{\VariableText}}

\newenvironment{titled}[1]
     {\def\VariableText{{#1}}\begin{VariableNoNum}}
     {\end{VariableNoNum}}

\newenvironment{NumberedSubSection}[1]
     {\def\VariableText{{\textrm{\textbf{#1}}}}\begin{VariableBold}}
     {\end{VariableBold}}

\def\Math#1{\def\MathString{#1}\futurelet\MathDelim\MathChoose}
\def\MathChoose{\ifmmode\let\MathDo\MathString%
              \else\let\MathDo\MathSkip\fi%
              \MathDo}
\def\MathSkip{\ifx\MathDelim/\def\MathDo{$\MathString$\EatOne}%
              \else\def\MathDo{$\MathString$}\fi%
              \MathDo}

\def\Text#1{\def\TextString{#1}\futurelet\TextDelim\TextSkip}
\def\TextSkip{\ifx\TextDelim/\def\TextDo{\TextString\EatOne}%
              \else\let\TextDo\TextString\fi%
              \TextDo}
\def\EatOne#1{}

\newdimen\tmpdimen
 \def\justboxit#1{{\setbox0\hbox{#1}\tmpdimen\dp0
   \advance\tmpdimen 0.4pt \lower\tmpdimen 
   \vbox{\hrule\hbox{\vrule\unhbox0\vrule}\hrule}}}
 \def\justaddborder#1#2{{\setbox0\hbox{#1}\tmpdimen\dp0
   \advance\tmpdimen #2\lower\tmpdimen
   \vbox{\kern#2\hbox{\kern#2\unhbox0\kern#2}\kern#2}}}
\def\CenterBox#1{\centerline{\justboxit{\justaddborder {#1}{3pt}}}}

\newcommand{\Ext}{\operatorname{Ext}}    
\newcommand{\colim}{\operatorname{colim}}
\newcommand{\Steen}{\Math{\mathcal A}}   
\newcommand{\Steenp}{\Math{\Steen_p}}    
\newcommand{\R}{\mathbb{R}}              
\newcommand{\Z}{\mathbb{Z}}              
\newcommand{\Fp}{\mathbb{F}_p}           
\newcommand{\Fq}{\mathbb{F}_q}           
\newcommand{\Ftwo}{\mathbb{F}_2}         
\newcommand{\tensor}{\otimes}            
\newcommand{\res}{\operatorname{res}}    
\newcommand{\tr}{\operatorname{tr}}      
\newcommand{\weq}{\simeq}                
\newcommand{\iso}{\cong}                 
\newcommand{\kP}{kP}                     
\newcommand{\FpP}{\Fp[P]}                
\newcommand{\RightArrow}[1]{\xrightarrow{#1}}    
\newcommand{\normal}{\lhd}               
\newcommand{\Syl}{\operatorname{Syl}}    
\newcommand{\cx}[2]{\operatorname{cx}_{#1}(#2)}  
\newcommand{\mi}{\mathfrak{m}}           
\newcommand{\Fam}{\Math{\mathcal F}}     
\newcommand{\pcgroup}{$p$C--group}       
\newcommand{\fib}{\vrule width 1.75cm height .1pt}

\newcommand{\Hom}{\operatorname{Hom}}

\begin{document}
\title[Cohomology of Finite Groups]{Lectures on the Cohomology of Finite Groups}

\author[Alejandro Adem]{Alejandro Adem$^*$}

\address{Department of Mathematics,
         University of British Columbia,
         Vancouver BC V6T 1Z2, Canada}

\email{adem@math.ubc.ca}
\thanks{{$^*$}Partially supported by NSF and NSERC. The author is
very grateful to Bill Dwyer for his generous help in
preparing this manuscript.}
\date{\today}
\abstract These are notes based on lectures given at the
summer school ``Interactions Between Homotopy
Theory and Algebra'', which was held at
the University of Chicago in the summer of 2004.
\endabstract
\maketitle \tableofcontents

\section{Introduction}

Finite groups can be studied as groups of symmetries in different
contexts. For example, they can be considered as groups of
permutations or as groups of matrices. In topology we like to think
of groups as transformations of interesting topological spaces,
which is a natural extension of the classical problem of describing
symmetries of geometric shapes. It turns out that in order to
undertake a systematic analysis of this, we must make use of the
tools of homological algebra and algebraic topology. The context for
this is the \textsl{cohomology of finite groups}, a subject which
straddles algebra and topology. Groups can be studied homologically
through their associated group algebras, and in turn this can be
connected to the geometry of certain topological spaces known as
\textsl{classifying spaces}. These spaces also play the role of
building blocks for stable homotopy theory and they are ubiquitous
in algebraic topology.

In these notes we have attempted to lay out a blueprint for the
study of finite group cohomology which totally intertwines the
algebraic and topological aspects of the subject. Although this may
pose some technical difficulties to the readers, the advantages
surely outweigh any drawbacks, as it allows them to understand the
geometric motivations behind many of the constructions.

The notes reflect the content of the lectures given by the author at
the summer school \textsl{Interactions Between Homotopy Theory and
Algebra} held at the University of Chicago in August 2004. The first
talk was essentially an outline of finite group cohomology from a
combined algebra/topology point of view; the second talk was a
discussion of connections between group cohomology and
representation theory via the concept of minimal resolutions;
finally the third talk was a discussion of the role played by group
cohomology in the study of transformation groups. This is reflected
in the organization of this paper; note that the final section
summarizes recent work by the author and others on constructing free
group actions on a product of two spheres. The style is somewhat
terse, as the goal was to expose the summer school participants to
research topics in the subject, beyond the usual background
material. However, there are several excellent books on the
cohomology of finite groups (listed in the bibliography), and the
hope is that the readers of this very short survey will feel
motivated to consult them, as well as the original research papers
in the literature.

The author is grateful to the organizers of the summer school for
being given the opportunity to give these lectures as well as for
their help in preparing this manuscript.

\section{Preliminaries}

We begin by recalling some basic facts about classifying spaces of
finite groups and group cohomology; useful references for this
material are \cite{Adem-Milgram}, \cite{Benson} and \cite{Brown}.

Let $G$ denote a finite group and $EG$ a free, contractible
$G$-complex (this means that $EG$ is a CW--complex on which $G$
acts in a way which freely permutes the cells). The orbit space
$BG=EG/G$ is called a \emph{classifying space} for $G$. The
quotient map $EG\to BG$ is a principal covering space with group
$G$, and from this it follows directly that $BG$ is a space of
type $K(G,1)$. In particular, $BG$ is path connected, the
fundamental group of $BG$ is isomorphic to $G$, and the higher
homotopy groups of $BG$ are trivial. The homotopy type of $BG$
does not depend upon the choice of $EG$.

\begin{rem}
  The term ``classifying space'' stems from the fact that isomorphism
  classes of principal $G$ bundles over a CW--complex $B$ correspond
  bijectively to homotopy classes of maps $B\to BG$ (see \cite{St}).
  In other words,
  such bundles are classified by maps to~$BG$. The correspondence
  sends a map $B\to BG$ to the pullback over $B$ of the projection
  $EG\to BG$.
\end{rem}

\begin{defn}
  For any $\Z[G]$--module $M$, we define the \emph{cohomology groups
  of $G$ with
  coefficients in $M$}, denoted $H^i(G,M)$ ($i\ge0$) by the formula
  \[
          H^i(G,M) = H^i(BG,M)\,.
  \]
\end{defn}

The notation ``$H^i(BG;M)$'' above signifies singular (or cellular) cohomology with
local coefficients; the local coefficient system is derived from the
given action of $G=\pi_1BG$ on~$M$.

\begin{rems}\label{BasicProperties}
  \par\noindent(1) Let $S_*(EG)$ be the integral singular complex of $EG$, so
  that essentially by definition, $H^i(G,M)$ is isomorphic to
  the $i$'th cohomology group of the cochain complex
  $\Hom_{\Z[G]}(S_*(EG),M)$. The fact that $EG$ is acyclic and that
    $G$ acts freely on $EG$ implies that the augmentation map
    $S_*(EG)\to\Z$ gives a free resolution of $\Z$ over~$\Z[G]$ (here
    $\Z$ serves as a trivial $G$-module, i.e., a module on which each
    element of $G$ acts as the identity). By basic homological
    algebra, we have
    \[
             H^i(G,M)=\Ext^i_{\Z[G]}(\Z, M)\quad\quad(i\ge 0)\,.
    \]
    In particular, $H^0(G,M)$ is the fixed submodule $M^G$, and the
    higher cohomology groups result from applying to $M$ the higher
    right derived functors of this fixed--point construction.

   \par\smallskip\noindent(2) Using joins, we may construct a model for $EG$
   which is functorial in $G$, namely $EG=\colim_iG^{*i}$, where $G^{*i}$
     is the join $G*G*\cdots*G$, $i$~times. The points of $EG$ can be
     thought of as infinite formal sums $\sum_{i\ge0}t_ig_i$, where
     $g_i\in G$, $t_i\in[0,1]$, only finitely many $t_i$ are nonzero,
     and $\sum t_i=1$. We define a right $G$-action on $EG$ by
    \[
         (\sum_it_ig_i)\cdot g =\sum_i t_i(g_ig)\,.
    \]
    From this point of view, the space $EG$ is endowed with the
    smallest topology which makes the coordinate functions
    \[
        t_i:EG\to [0,1] \quad \text{ and } \quad g_i:t_i^{-1}(0,1]\to G
    \]
    continuous. Taking joins increases connectivity, so $EG$ is
    contractible. It is clear that $G$ acts freely on $EG$, and it is
    an interesting exercise to construct a CW structure on $EG$ so
    that the action of $G$ permutes the cells. If $f:H\to G$ is a group
     homomorphism, we get induced maps $EH\to EG$ and $BH\to BG$, as
     well as maps $H^i(G,M)\to H^i(H,M|_H)$. Here $M|_H$ is~$M$, treated as an
     $\Z[H]$-module by means of the ring homomorphism
     $\Z[f]:\Z[H]\to\Z[G]$.

   \par\smallskip\noindent(3) If $R$ is a ring with a trivial action
   of $G$, then $H^*(G,R)$ has a natural
   graded multiplicative structure, arising from the cup product.

  \par\smallskip\noindent(4) If $\Fp$ is given the trivial $G$--module
  structure, then $H^*(G,\Fp)$ is the ordinary mod~$p$ cohomology ring
  of $BG$, and so it has an action of the Steenrod algebra \Steenp.

  \par\smallskip\noindent(5) The low dimensional cohomology groups $H^i(G;M)$,
  $i=1,2$ arose classically in the study of extensions of $G$ by
  $M$. Such an extension is a short exact sequence
  \[
         1 \to M \to \tilde G \to G \to 1
  \]
  in which the conjugation action of $\tilde G$ on $M$ induces the
  given $G$-module structure on $M$. Isomorphism classes
  of these extensions correspond bijectively to $H^2(G,M)$. If the
  extension corresponds to the zero cohomology class, then there is a
  section $G\to \tilde G$, and $H^1(G,M)$ acts freely and transitively
  on the set of $M$-conjugacy classes of such sections.

  The group $H^3(G;M)$ arises in a similar but more complicated
  context in studying extensions of $G$ by a nonabelian group $H$; the
  relevant $M$ is the center of $H$, and elements of $H^3(G,M)$ come
  up as obstructions to realizing a map from $G$ to the outer
  automorphism group of $H$ via an extension.
\end{rems}

\begin{exs} (1) If $G=\Z/2$, then $BG$ is equivalent to $\R P^\infty$,
  and
  $H^*(BG,\Ftwo)$ is a polynomial algebra $\Ftwo[x]$ with
  $|x|=1$.

\par\smallskip\noindent(2) If $G=\Z/p$ with $p$~odd, then $BG$ is
equivalent to the infinite lens space
$L_p^\infty=\colim_iS^{2i-1}/\mu_p$, and $H^*(BG,\Fp)$ is the tensor
product $\Lambda(x)\tensor\Fp[y]$ of an exterior algebra and a
polynomial algebra, where $|x|=1$ and $|y|=2$. In this case $y$ is
the Bockstein $\beta(x)$.

\par\smallskip\noindent(3) It is easy to see that $E(G\times H)$ can
be taken to be $EG\times EH$, so that $B(G\times H)$ is homotopy
equivalent to $BG\times BH$. By the K\"unneth formula, then, there are
isomorphisms
\[
H^*((\Z/p)^n,\Fp)\cong \begin{cases}

                        \Fp[x_1,\ldots,x_n] & p=2\\
                    \Lambda(x_1,\ldots,x_n)\tensor\Fp[y_,\ldots,y_n]
                    &p\text{ odd.}
                   \end{cases}
\]
\end{exs}

\begin{NumberedSubSection}{Restrictions and transfers}
Let $H\subset G$ be a subgroup, and note that $EG$ is a contractible
space with a free, cellular action of $H$, so that $BH\simeq EG/H$.
Hence (cf. \ref{BasicProperties}) we have a map
\[
     BH = EG/H \to EG/G = BG\,.
\]
This induces a \emph{restriction map} $\res^G_H:H^*(G,M)\to
H^*(H,M|_H)$.
Note that if $R$ is a ring with a trivial action of $G$, then the
restriction map is a map of graded $R$-algebras.

Now if we consider the cell structure of $EG$, then over any cell
$\sigma$ of $BG$ there are $[G:H]=\#(G/H)=n$ cells
$g_1\tilde\sigma,\ldots,g_n\tilde\sigma$ of $BH$, where $\tilde\sigma$
is some fixed chosen cell above $\sigma$ and $g_1,\ldots,g_n$ are
coset representatives for $G/H$. We define
\[
   \psi:C_*(BG,\Z) \to C_*(BH,\Z)
\]
by setting $\psi(\sigma)=\sum_{i=1}^ng_i\tilde\sigma$. This idea can
be exploited in slightly greater generality to construct a transfer
map $\tr^G_H:H^*(H,M|_H)\to H^*(G,M)$.

There are a few basic formulas involving these maps.
\begin{enumerate}
\item If $H\subset K\subset G$, then $\res^K_H\cdot\res^G_K=\res^G_H$.
\item $\tr^G_H\cdot\res^G_H(x)=[G:H]x,\quad\forall x\in H^*(G,M)\,.$
\item Suppose that $H$ and $K$ are subgroups of $G$, and that $G$ is
  written as a disjoint union $\cup_{i=1}^m Hg_iK$ of double cosets.
  Then
  \[ \res^G_H\cdot\tr^G_K=\sum_{i=1}^m\tr^H_{H\cap
    K^{g_i}}\cdot\res^{K^{g_i}}_{H\cap K^{g_i}}\cdot C_{g_i}\,,
  \]
  where $C_{g_i}$ is induced by conjugation with~$g_i$.
\end{enumerate}

\begin{exercise}
  Let $P$ be a Sylow $p$-subgroup of $G$, and $N_G(P)$ its normalizer
  in~$G$. If $P$ is abelian, show that restriction induces an
  isomorphism
  \[  H^*(G,\Fp)\cong H^*(N_G(P),\Fp)\,.
   \]
\end{exercise}

\end{NumberedSubSection}

\begin{NumberedSubSection}{Lyndon--Hochschild--Serre spectral
    sequence}
If $H$ is a normal subgroup of $G$, then $G/H$ acts freely on
$BH\weq EG/H$, and so we have a fibration
\[
\begin{CD}
      BH @>>> E(G/H) \times_{G/H} EG/H \weq BG\\
         @.                      @VVV\\
                  @.            B(G/H)
\end{CD}
\]
This gives a spectral sequence
\[
 E_2^{p,q} = H^p(G/H, H^q(H,M))\Rightarrow H^{p+q}(G,M)\,,
\]
in which the $E_2$-page involves local coefficient cohomology
determined by an action of $G/H$ on $H^*(H,M)$.

\begin{exs}
  (1) Let $G$ be the alternating group $A_4$. There is a group
extension
\[ 1\to (\Z/2)^2 \to G=A_4\to \Z/3 \to 1\,.
\]
The associated mod~2 cohomology spectral sequence collapses to give
the formula
\[
\begin{aligned}
  H^*(A_4,\Z/2)&\iso\Ftwo[x_1,x_3]^{\Z/3}\\
               &\iso \Ftwo[u_2,v_3,w_3]/(u_2^3+v_3^2+w_3^2+v_3w_3)\,.
\end{aligned}
\]

\par\smallskip\noindent(2) Let $G$ be the dihedral group $D_8$, which
can be written as the wreath product $\Z/2\wr\Z/2$ or equivalently as
the semidirect product $(\Z/2\times\Z/2)\rtimes\Z/2$. The mod~2
cohomology spectral sequence of the semidirect product extension
collapses strongly at $E_2$ and gives the formula
\[
H^*(D_8,\Ftwo)\iso \Ftwo[x_1, e_1, y_2]/(x_1e_1)\,.
\]
\end{exs}

Given any $G$, we can find a monomorphism $G\to U(n)$ to obtain a
fibration
\[
        U(n)/G \to BG \to BU(n)\,.
\]
Here $U(n)$ is the unitary group of rank~$n$, and we are
implicitly referring to the fact that classifying spaces can be
constructed not just for finite groups but also for topological
groups such as $U(n)$. Recall that $H^*(BU(n),\Fp)$ is isomorphic
to a polynomial algebra $\Fp[c_1,\ldots,c_n]$ (the generators are
the universal Chern classes). Venkov \cite{Venkov} used the above
fibration sequence to prove
\begin{thm}
If $G$ is a finite subgroup of $U(n)$, then $H^*(G,\Fp)$ is a finitely
generated module over $H^*(BU(n),\Fp)$, and its Poincar\'e series is a
rational function of the form
\[
p_G(t)=
\sum\dim H^i(G,\Fp)t^i=\frac{r(t)}{\prod_{i=1}^m(1-t^{2i})}\,,
\]
where $r(t)\in\Z[t]$.
\end{thm}

\begin{ex}
  Let $G=Q_8$, the quaternion group of order~$8$. Then $G\subset
  SU(2)$, and we have a fibration
  \[
     SU(2)/G \to BQ_8 \to BSU(2)\,.
  \]
Here $H^*(BSU(2),\Ftwo)\iso \Ftwo[u_4]$ and $H^*(BQ_8,\mathbb F_2)$
is a free
$\Ftwo[u_4]$-module, with basis given by
\[H^*(SU(2)/Q_8)\cong
\Ftwo[x_1,y_1]/(x_1^2+x_1y_1+y_1^2, x_1^2y_1+x_1y_1^2)\,.\]
Then $H^*(Q_8)/(u_4)\iso H^*(SU(2)/Q_8, \mathbb F_2)$.
\end{ex}

\begin{titled}{Question}
What is the order of the pole of $p_G(t)$ at $t=1$? This is
known as the Krull
dimension of $H^*(G,\Fp)$.
\end{titled}

\begin{titled}{Answer}
The order of the pole of $p_G(t)$ at $t=1$ is the \emph{$p$-rank} $r_p(G)$ of
$G$, defined as
\[
  r_p(G)=\max\{n\mid (\Z/p)^n\subset G\}\,.
\]
\end{titled}

In the next section, we will try to explain this answer using
representation theory.
\end{NumberedSubSection}
\section{Minimal resolutions}

Let $P$ be a finite $p$-group, $k$ a field of characteristic $p$, and
$M$ a finitely generated $\kP$-module. Using the nilpotence of the
augmentation ideal $I\subset \kP$, one can see that the rank of the
projective cover of $M$ is $r=\dim M/IM =\dim H_0(P,M)$. Hence we have
an exact sequence
\[
   0 \to \Omega^1(M) \to (\kP)^r \to M \to 0
\]
(this sequence defines $\Omega^1(M)$) where the right hand map induces
an isomorphism $H_0(P,(\kP)^r)\iso (\Fp)^r\to H_0(P,M)$.
The long exact homology sequence associated to the above short exact
sequence shows that $\dim H_0(P,\Omega^1(M))=\dim H_1(P,M)$; this last
number is then in turn the rank of the projective cover of
$\Omega^1(M)$. Continuing on like this, and identifying
\[
           H_i(P,M)\iso H^i(P,M^*) \quad\quad\text{where }
           M^*=\Hom(M,k)\,,
\]
we obtain the following statement.

\begin{prop}
  If $P_*\to M$ is a minimal projective resolution for $M$, then
\[
   \dim_{\Fp} P_i = \vert P\vert \dim H^i(P,M^*)\,.
\]
\end{prop}

\begin{cor}
  The following are equivalent:
  \begin{enumerate}
  \item $M$ is projective,
  \item $H^r(P,M)=0$ for some $r>0$, and
  \item $H^i(P,M)=0$ for all $i>0$.
  \end{enumerate}
\end{cor}

We want to determine projectivity by restriction to certain
subgroups. We will use group cohomology to do this, by applying
the following basic result of Serre (\cite{Serre}):

\begin{thm}
  Let $P$ denote a finite $p$-group which is not elementary
  abelian. Then there exist non--zero elements
  $\gamma_1,\ldots,\gamma_n\in H^1(P,\Fp)$ such that 
$$\beta(\gamma_1)\beta(
\gamma_2)\cdots\beta(\gamma_n)=0,$$ 
where $\beta$ is the Bockstein.
\end{thm}

\noindent Now if $\gamma\in H^1(\Z/p,\Z/p) $ corresponds to the identity
homomorphism $\Z/p\to\Z/p$, then its Bockstein in
$H^2(\Z/p,\Z/p)\iso\Ext^2_{\Fp[\Z/p]}(\Z/p,\Z/p)$ is the sequence
\[
  \Fp\to \Fp[\Z/p]\RightArrow{t-1}\Fp[\Z/p]\to\Fp\,.
\]
For non--zero $\gamma_i\in H^1(P,\Fp)$, we take a homomorphism $\phi:P\to\Z/p$
representing it, with kernel $H_i\normal P$ (a maximal subgroup).
Pulling back the representative for the Bockstein of the identity map
of $\Z/p$ gives an expression for $\beta\gamma_i$ as a sequence
\[
  \Fp\to\Fp[P/H_i]\to \Fp[P/H_i] \to \Fp\,.
\]
Taking products of cohomology classes corresponds to splicing
extensions, so that $\beta(\gamma_1)\cdots\beta(\gamma_n)$ is
represented by the sequence
\[
\begin{aligned}
  \Fp\to \Fp[P/H_n]\to \Fp[P/H_n] &\to\Fp[P/H_{n-1}]
\to \Fp[P/H_{n-1}]\to \\
\cdots&\to
\Fp[P/H_1]
\to
\Fp[P/H_1]\to \Fp\,.
\end{aligned}
\]
We leave off the copies of $\Fp$ on either end and interpret this as a
cochain complex $C^*$, with $H^*(C^*)=H^*(S^{2n-1},\Fp)$.

\begin{prop}
  If $P$ is not elementary abelian, then the $\FpP$--module  $M$ is
  projective if and only if $M\vert_K$ is projective for each maximal
  subgroup $K\subset P$.
\end{prop}

\begin{proof}
  Consider $C^*\otimes M$; we will compute $H^*(P, C^*\otimes M)$ in
  two different ways (what we are computing is sometimes called the
  \emph{hypercohomology of $P$ with coefficients in $C^*\otimes
    M$}). There are two spectral sequences
  \begin{enumerate}
  \item $E_2^{r,q}=H^r(P, H^q(C^*\otimes M))\Rightarrow
    H^{r+q}(P,C^*\otimes M)$, and
   \item $E_1^{r,q}=H^q(P, C^r\otimes M)\Rightarrow H^{r+q}(P,
     C^*\otimes M)$.
  \end{enumerate}
  The {first} one has a single differential determined by
  $\beta(\gamma_1)\cdots\beta(\gamma_n)=0$ (by construction), and so
  it collapses at $E_2$ giving
  \[
     H^\ell(P,C^*\otimes M)\iso H^\ell(P,M)\oplus H^{\ell-2n+1}(P,M).
\]
  For the second spectral sequence, note that
  \[
      H^q(P,C^r\otimes M)\iso H^q(H_i,M)\,,
   \]
  where $H_i\subset P$ is maximal. So if $M\vert_K$ is projective for
  any maximal $K\subset P$, we get that $H^t(P,C^*\otimes M)$ vanishes
  for $t>\!\!>0$. Combining these two calculations shows that
  $H^t(P,M)=0$ for $t>\!\!>0$, and so $M$ is projective. The opposite
  implication is clear.
\end{proof}

An immediate consequence is

\begin{thm}
  (Chouinard) An $\Fp[G]$--module $M$ is projective if and only if
  $M\vert_E$ is projective for all elementary abelian $p$-subgroups of~$G$.
\end{thm}

\begin{proof}
  Let $P\in\Syl_p(G)$; if $M\vert_P$ is free (equivalently,
  projective), so is the module
  $\Fp[G]\otimes_{\Fp[P]}M$; as this induced
  module contains $M$ as a direct summand (use the fact that the index
  of $P$ in $G$ is prime to~$p$), it follows that $M$ is
  projective. This shows that $M$ is projective if and only if $M\vert_P$
  is projective; now apply our previous result repeatedly to reduce
  the problem of testing projectivity for $M$ to the problem of
  examining $M\vert_E$ for every elementary abelian subgroup of~$P$.
\end{proof}

If $V_*$ is a graded $k$-vector space, we can define its \emph{growth
  rate} $\gamma(V_*)$ by
\[
  \gamma(V_*)=\min\left\{n\ge 0\mid \lim_{t\to\infty}\frac{\dim
    V_t}{t^n}=0\right\}.
\]

\begin{defn}
The \emph{complexity} $\cx GM$ of a $kG$-module $M$ is defined to be
$\gamma(k\otimes_{kG}P_*)$, where $P_*$ is a minimal projective
resolution of $M$ over $kG$.
\end{defn}

\begin{thm}
  (Quillen/Alperin--Evens \cite{Q1}, \cite{Alperin-Evens})
  \[
    \cx GM = \max_{E\subset G}\left\{\cx E{M\vert_E}\right\}\,,
  \]
  where $E$ runs through elementary abelian $p$-subgroups of~$M$.
\end{thm}

We can sketch a proof of this: as before we can reduce it to
$p$-groups. Now for a $p$-group $P$ we have that
\[
\begin{aligned}
  \cx PM &= \gamma(H^*(P,M^*))\\
         &= \max_E \gamma(H^*(E,M^*))\\
         &= \max_E \cx EM\,.
\end{aligned}
\]
Going from the first to second line here uses our previous argument
and Serre's result. \qed

For trivial coefficients, this implies that the Krull dimension of
$H^*(G,\Fp)$ is precisely $r_p(G)$, as we have an explicit computation
of $H^*(E,\Fp)$.

Let $V_G(k)$ denote the maximal ideal spectrum for $H^*(G,k)$. The
restriction $\res^G_E$ induces a map
\[
(\res^G_E)^* : V_E(k)\to V_G(k)\,.
\]

\begin{thm}
  (Quillen \cite{Q1}) Let $A_p(G)$ be the set of all elementary abelian
  $p$-subgroups of $G$. Then
\[
V_G(k)=\cup_{E\in A_p(G)} (\res^G_E)^* (V_E(k))\,.
\]
\end{thm}

We can view a maximal ideal $\mi\in V_G(k)$ as the kernel of a nonzero
homomorphism $H^*(G,k)\RightArrow\alpha\bar k$, where $\bar k$ is an
algebraic closure of $k$. Quillen's theorem says that every such
homomorphism $\alpha$ is of the form
\[H^*(G,k)\RightArrow{\res^G_E} H^*(E,k) \RightArrow\beta \bar k\]
for some $E\in A_p(G)$.

For more details on the methods outlined in this section, we
refer the reader to the paper by J. Carlson \cite{Carlson}.

\section{Computations and further structure}

Let \Fam/ denote a family of subgroups of $G$, i.e., if $H\in\Fam$,
$H'\subset H$, then $H'\in\Fam$ and $gHg^{-1}\in\Fam$ for any $g\in
G$.
Then we can define
\[
\lim_{H\in\Fam} H^*(H,\Fp)=\left\{(\alpha_H)\text{ such that }
\begin{aligned} \alpha_{H'}&=\res^H_{H'}\alpha_H\text{ if } H'\subset
H\\
\alpha_{H'}&=c_g\alpha_H\text{ if } H'=gHg^{-1}
\end{aligned}
\right\}
\]
We can use this construction to reinterpret our previous results.

\begin{thm}
  (Cartan--Eilenberg \cite{Cartan-Eilenberg})
  Let $S_p(G)$ denote the family of all
  $p$-subgroups of $G$. Then the restrictions induce an
  \emph{isomorphism}
\[ H^*(G,\Fp)\iso\lim_{P\in S_p(G)} H^*(P,\Fp)\,.\]
\end{thm}

\begin{thm}
  (Quillen--Venkov \cite{Quillen-Venkov})
  Let $A_p(G)$ denote the family of all
  $p$-elementary abelian subgroups of $G$. Then the restrictions
  induce an $F$-isomorphism\footnote{By which we mean that the
  kernel of $\theta$ is nilpotent and that a sufficiently high power of any
  element in the target lies in the image of $\theta$.}
\[
  \theta: H^*(G,\Fp)\to\lim_{E\in A_p(G)} H^*(E,\Fp)\,.\]
\end{thm}

We need to compute $H^*(G,\Fp)$ for interesting classes of groups. In
many cases this will involve finding a collection of subgroups
$H_1,\ldots,H_\ell$ such that the map
\[
H^*(G)\to \oplus_{i=1}^\ell H^*(H_i)
\]
is injective, in which case this collection is said to \emph{detect}
the cohomology.

\begin{NumberedSubSection}{Calculational methods}
Calculating the cohomology of finite groups can be quite challenging,
as it will involve a number of complicated ingredients. We outline the
main techniques and provide some examples.

How to compute:

\par\smallskip\noindent Step 1. Reduce to the Sylow $p$-subgroup via
the Cartan--Eilenberg result, and then combine information about the
cohomology of $p$-groups with stability conditions.

\par\smallskip\noindent Step 2. Determine $A_p(G)$ and use Quillen's
result to compute
$$H^*(G)/\sqrt{O},$$ where $\sqrt{O}$ is the
radical
    of the ideal of elements of strictly positive degree.
This will require computing rings of invariants of the form
$H^*(E)^N$, where $E$ is $p$-elementary abelian and $N$ is the
normalizer of $E$.

\par\smallskip\noindent Step 3. Let $\vert A_p(G)\vert$ be the $G$-CW
complex which realizes the poset $A_p(G)$. A result due to Webb
\cite{Webb} implies that $H^*(G,\Fp)$ can be computed from the
cohomology of the groups  $N_G(E)$ and the cohomology of certain
of their intersections. These are the intersections which appear
as simplex stabilizers for the action of $G$ on $|A_p(G)|$.) More
precisely, the mod $p$ Leray spectral sequence associated to the
projection $EG\times_G|A_p(G)|\to |A_p(G)|/G$ collapses at $E_2$
with
\[
   E_2^{s,t}=\begin{cases}
                 H^t(G) & s=0\\
                   0 & s>0
              \end{cases}
\]
If $G$ is of Lie type, we can use the Tits building and information
about the parabolic subgroups of $G$ to facilitate this calculation,
as it will be equivariantly equivalent to the poset space above.

\begin{ex}
  Let $S_n$ be the symmetric group of degree~$n$. In this case we have
  that Quillen's map induces a mod~2 isomorphism
  \[
  H^*(S_n,\Ftwo)\RightArrow{\iso} \lim_{E\in A_2(S_n)} H^*(E,\Ftwo)\,.
\]
This means that the cohomology is detected on elementary abelian
subgroups. To determine $H^*(S_n)$ we need to glue together the
different bits from these detecting subgroups. This can be done
for a general $S_n$, using invariant theory and combinatorics. In
fact, according to Milgram \cite{Adem-Milgram} and Feshbach
\cite{Feshbach}, $H^*(S_\infty)$ surjects onto $H^*(S_n)$, and
$H^*(S_\infty)$ is known (from Nakaoka \cite{Nakaoka}). As an
explicit example, $H^*(S_4)\iso \Ftwo[x_1,y_2,c_3]/(x_1c_3)$, is
detected by a map to
\[
H^*(S_4) \to
      H^*(V_1)^{\Z/2}\iso \Ftwo[\sigma_1,\sigma_2]
\quad\oplus\quad
     H^*(V_2)^{GL_2(\Ftwo)}\iso \Ftwo[d_2,d_3]
\]
where
\[
   \begin{aligned}
   x_1&\mapsto (\sigma_1,0)\\
   y_2&\mapsto (\sigma_2,d_2)\\
   c_3&\mapsto (0,d_3)
   \end{aligned}
\]
\end{ex}

\begin{rems}
  \par\noindent(1) $H^*(A_n,\Ftwo)$ can be obtained from the cohomology
  of $H^*(S_n,\Ftwo)$.

  \par\smallskip\noindent(2) If $G$ is one of the standard algebraic groups
  over~$\Z$,  then $H^*(G(\Fq),\Fp)$, $p\ne q$, was computed by
  Quillen \cite{Q2} and others. In these cases, the cohomology is
  \emph{detected} on abelian subgroups.

\par\smallskip\noindent(3) Quillen \cite{Q3}
also computed the mod~2 cohomology of the
extra--special $2$-groups, described as central extensions
\[
   1\to\Z/2\to P\to (\Z/2)^r\to 1\,.
\]
This can be done with the Eilenberg-Moore spectral sequence.
\end{rems}

We now give some examples of calculations for simple groups at $p=2$.

\begin{exs}
  \par\noindent(1) $G=L_3(2)$: we use the Tits
building, which is equivalent to $|A_2(G)|$;
in this case it has quotient

$$\xymatrix{
&S_4~~{\buildrel{D_8}
\over{\bullet\!\!\hbox{\fib}\!\!\bullet}}~~S_4\\}
$$

We compute the cohomology of $G$
as the intersection of two
copies of the cohomology of $S_4$ in the cohomology
of $D_8$:
\[
       H^*(L_3(2)) \iso \Ftwo[u_2,v_3,w_3]/(v_3w_3)\,.
\]
\par\smallskip\noindent(2) $G=A_6$: here we use $H^*(S_6)/(\sigma_1)=H^*(A_6)$,
yielding
\[
     H^*(A_6) = \Ftwo[\sigma_2,\sigma_3,c_3]/(c_3\sigma_3)\,.
\]
\par\smallskip\noindent(3) $G=M_{11}$, the first Mathieu group;
in this case
  $|A_2(G)|/G$ looks like:

$$\xymatrix{
&S_4~~{\buildrel{D_8}
\over{\bullet\!\!\hbox{\fib}\!\!\bullet}}~~GL_2(3)\\}
$$

\noindent giving
\[
     H^*(M_{11}) \iso \Ftwo[v_3,u_4, w_5]/(w_5^2+v_3^2u_4)\,.
\]
\par\smallskip\noindent(4) (Adem--Milgram \cite{Adem-Milgram2})
If $S$ is a simple group of rank~$3$ or less at $p=2$, then
$H^*(S,\Ftwo)$ is Cohen--Macaulay. This means that the cohomology
ring is a finitely generated free module over a polynomial
subring.
\end{exs}
\end{NumberedSubSection}

From the calculation above, we see that the mod 2 cohomologies
of $A_6$ and $L_3(2)$ are isomorphic, even though there is no
non--trivial homomorphism between these two simple groups.
However, there is an infinite amalgam
$\Gamma=\Sigma_4*_{D_8}\Sigma_4$ such that
there are homomorphisms $\Gamma\to A_6$ and $\Gamma\to L_3(2)$
inducing mod 2 homology equivalences. From this one can
show that the classifying
spaces $BA_6$ and $BL_3(2)$ are equivalent at $p=2$.

\begin{NumberedSubSection}{Depth and detection}

An important invariant for the mod~$p$ cohomology $H^*(G,\Fp)$ is
its \emph{depth}, defined as the maximal length of a regular
sequence. We have the following basic result of Duflot
\cite{Duflot}.

\begin{thm}
  If $Z(G)$ is the center of $G$ and $r_p(Z(G))=z$, then
  \[
       \operatorname{depth} H^*(G) \ge z\,.
\]
\end{thm}

We will sketch a proof of this, for a $p$-group~$P$.
Let $Z(P)\iso \Z/p^{n_1}\times\cdots\times\Z/p^{n_z}$, and for
each summand choose a $1$-dimensional faithful complex
representation~$\chi_i$; this extends to a representation of
$Z(P)$. Let $V_i=\operatorname{Ind}_{Z(P)}^P(\chi_i)$ and consider
its associated sphere $S(V_i)$. Then, if
$X=S(V_1)\times\cdots\times S(V_z)$, $X$ has an action of~$P$.
Consider the mod $p$ Serre spectral sequence for the fibration
\[
\begin{CD}
  X @>>> EP\times_P X\\
  @.            @VVV\\
     @.     BP
\end{CD}
\]
Let $s=2[P:Z]$.  Then $H^*(X)\iso\Lambda_{\Fp}(u_1,\ldots,u_z)$, where
the degree of $u_i$ is $s-1$; each class $u_i$ is invariant under the
action of $P$ on $H^*(X)$ and so represents a class in $E_2^{0,s-1}$.
For positional reasons $u_i$ survives to $E_s^{0,s-1}$, and we write
$d^s(u_i)=\alpha_i\in H^{s}(P)$. The key result which can be proved
inductively is that $\alpha_1,\ldots,\alpha_z$ form a \emph{regular
  sequence}. This implies that there is an isomorphism
\[
H^*(EP\times_PX)\iso H^*(P)/(\alpha_1,\ldots,\alpha_z)\,.
\]
There is a special case: if $r_p(Z(P))=r_p(P)$, then every element of
order~$p$ in $P$ is central. In this case
\begin{itemize}
\item $H^*(P)$ is Cohen-Macaulay, being free and finitely generated
  over $\Fp[\alpha_1,\ldots,\alpha_z]$.
\item $P$ acts freely on $X$ and there is an isomorphism \[H^*(X/P)\iso
  H^*(P)/(\alpha_1,\ldots,\alpha_z)\,.\]
\end{itemize}
Note that we have a special geometric orientation class
\[
\mu\in H^{\text{top}}(X/P)\,.
\]
This element can be pulled back and used to construct an
``undetectable'' class in $H^*(P)$, yielding

\begin{thm}
  (Adem--Karagueuzian \cite{AK}) If every element of order~$p$ in $P$ is central
  (i.e., $P$ is a \pcgroup),
  then $H^*(P)$ is Cohen--Macaulay and $H^*(P)$ is undetectable:
there exists an element $x\ne 0$ in $H^*(P)$ such that
$res^P_H(x)=0$ for all proper subgroups $H\subset P$.

\end{thm}

On the other hand, we have

\begin{thm}
  (Carlson \cite{Carlson2}) If $H^*(G,\Fp)$ has depth $r$, then this cohomology is
  detected by subgroups of the form $C_G(E)$, where $E\subset G$ is
  $p$-elementary abelian of rank~$r$.
\end{thm}

\begin{rems}
  In particular, this tells us that if $H^*(G,\Fp)$ is Cohen-Macaulay,
  then it is detected on subgroups of the form $C_G(E)$, where $E$ has
  maximal rank. These are necessarily \pcgroup{}s. Hence we have the
  converse statement, that if $H^*(G,\Fp)$ is Cohen--Macaulay and
  undetectable, then $G$ is a \pcgroup. We can think of \pcgroup s as
  ``universal detectors'' for Cohen-Macaulay cohomology rings. Hence
  determining their cohomology is a basic problem.
\end{rems}

\begin{ex}
  Consider a ``universal'' or maximal central extension
\[
  1 \to H_2E \to P \to E\to 1
\]
where $E=(\Z/p)^n$, $p$ is an odd prime,
and the differential $H^1(H_2(E))\to H^2(E)$ is an
isomorphism. Then one can show that if $p> {n\choose 2} +1$, the Serre
spectral sequence collapses at $E_3$ and we have an exact sequence
\[
0 \to \left(\zeta_1,\ldots,\zeta_{{n+1}\choose {2}}\right) \to H^*(P)\to
\operatorname{Tor}_{\Fp[c_{ij}]}(\Lambda(e_1,\ldots,e_n),\Fp)\to 0
\]
where the Tor term is determined by $c_{ij}\mapsto e_ie_j$ for
$i<j$, $i,j=1,\ldots,n$ (see \cite{AP}).
\end{ex}

\begin{problem}
  Motivated by the above example, we can raise the following
  question. Let
\[ 1 \to V\to P\to W\to 1
\]
 be a \emph{central} extension, where $V$ and $W$ are elementary
 abelian groups. Can the Eilenberg-Moore spectral sequence fail to
 collapse at
\[
   E_2=\operatorname{Tor}_{H^*(K(V,2))}(H^*(W),\Fp) \,\text{?}
\]
If so, give reasonable conditions on the $k$--invariants which imply a
collapse.
\end{problem}
\end{NumberedSubSection}

\begin{NumberedSubSection}{Duality for group cohomology}

We briefly recall an important condition on $H^*(G,\Fp)$, related to
duality. Let $k[\zeta_1,\ldots,\zeta_r]\subset H^*(G,k)$ be a
homogeneous system of parameters having
$\deg\zeta_i=n_i\ge2$.

\begin{thm}
  (Benson--Carlson \cite{Benson-Carlson})
  There exists a finite complex $C$ of projective
  $kG$--modules with $H^*(\Hom_k(C,k))\iso
  \Lambda(\bar\zeta_1,\ldots,\bar\zeta_r)$, with
  $\deg\bar\zeta_i=n_i-1$. There is a spectral sequence with
  \[
         E_2^{*,*} =
         H^*(G,k)\otimes\Lambda(\bar\zeta_1,\ldots,\bar\zeta_r)
  \]
  converging to $H^*(\Hom_{kG}(C,k)$, which satisfies Poincar\'e
  duality in formal dimension $s=\sum_{i-1}^r (n_i-1)$. In this
  spectral sequence we have $d_{n_i}(\bar\zeta_i)=\zeta_i$, and if
  $H^*(G,k)$ is Cohen--Macaulay, then
  \[
      H^*(G,k)/(\zeta_1,\ldots,\zeta_r)\iso H^*(\Hom_{kG}(C,k))\,.
  \]
\end{thm}

More succinctly,

\begin{thm}
Suppose that $G$ is a finite group. If the ring $H^*(G,k)$ is
Cohen-Macaulay, then it is Gorenstein with $a$-invariant zero.
\end{thm}

\begin{rems}
  In the Cohen--Macaulay case, $H^*(G)/(\zeta_1,\ldots,\zeta_r)$
  satisfies Poincar\'e duality, and its top degree is
  $\sum_{i=1}^r(n_i-1)=d$. The ``$a$-invariant'' is computed as
  $a=d-\sum_{i=1}^r(|\zeta_i|-1) =0$ in this case. Note the functional
  equation: if $p_G(t)$ is the Poincar\'e series for $H^*(G)$, then
  \[
        p_G(1/t) = (-t)^r p_G(t)\,.
  \]
  The previous theorem asserts that if $G$ is a finite group of
  rank~$r$, then there exists a projective $\Z G$ chain complex $C$
  with
  \[ H^*C^*\iso H^*(S^{n_1-1}\times\cdots\times S^{n_r-1})\,. \]
\end{rems}

In the next section we will describe the analogous problem in a
geometric setting, i.e., using free group actions on products of
spheres. This turns out to be much more difficult, as we shall see.
\end{NumberedSubSection}

\section{Cohomology and actions of finite groups}
We start by recalling a basic result.

\begin{thm}
  (P. A. Smith, 1940 \cite{Smith})
  If a finite group $G$ acts freely on a sphere,
  then all of its abelian subgroups are cyclic.
\end{thm}

What this means is that $G$ does not contain any subgroup of the form
$\Z/p\times\Z/p$ (``the $p^2$ condition'').

\begin{exs}
\par\noindent
  (1) $\Z/n$ acts freely on any $S^{2k+1}$.
\par\smallskip\noindent (2) The quaternion group
$Q_8$ acts freely on $S^3$.
\end{exs}

Later, Milnor \cite{Milnor} proved

\begin{thm}
  If a finite group $G$ acts freely on $S^n$. then every involution in
  $G$ is central.
\end{thm}

For example, the dihedral group $D_{2p}$ cannot act freely on any
sphere (``the $2p$ condition''). However, it turns out that this
is really a \emph{geometric} condition, as we have

\begin{thm}
  (Swan \cite{Swan}) A finite group $G$ acts freely on a finite complex
  $X\simeq S^n$ if and only if every abelian subgroup of $G$ is cyclic.
\end{thm}

Finally, using surgery theory, it was shown that

\begin{thm}
  (Madsen--Thomas--Wall \cite{MTW})
  $G$ acts freely on some sphere if and
  only if $G$ satisfies both the $p^2$ and the $2p$ condition for all
  primes~$p$.
\end{thm}

If $G$ acts freely on $X=S^n$, preserving orientation,
we have a Gysin sequence (with $\mathbb Z$ or $\mathbb F_p$
coefficients):

\[
\cdots\to H^{i+n}(X/G)\to H^i(G)\RightArrow{\cup x} H^{i+n+1}(G)\to
H^{i+n+1}(X/G)\to\cdots
\]
where $x\in H^{n+1}(G)$ is the Euler class of the action.
Note that $H^s(X/G)=0$ for $s>n$; this implies that
the map
\[
  \cup x: H^i(G)\to H^{i+n+1}(G)
\]
 is an isomorphism for all $i>0$. This in turn implies that $G$ has
 periodic cohomology, i.e., $H^*(G,\Fp)$ has Krull dimension one for
 all $p\mid |G|$.  It follows that $G$ satisfies the $p^2$-condition
 for all primes $p$.

 \begin{rems}
   (1) In fact, $G$ has periodic cohomology if and only if every
   abelian subgroup in $G$ is cyclic (Artin--Tate \cite{AT}).
\par\noindent(2) Given $G$ with periodic cohomology, say of minimal
period $d$, acting freely on $X\simeq S^n$, then $d\mid n+1$, but $d$
is not necessarily equal to $n+1$.
 \end{rems}

 \begin{ex}
   Let $S=Z/3\rtimes Q_{16}$, where the element of order~$8$ in $Q_{16}$
   acts nontrivially on $\Z/3$. Then $S$ has period \emph{four} but
   does not act freely on any finite homotopy $3$-sphere. Hence there
   is no closed $3$--manifold $M$ with $\pi_1(M)\iso S$.
 \end{ex}

 \begin{exs}
   \emph{Linear spheres.} Let $V$ be a unitary or orthogonal
   representation of $G$ which is fixed--point--free. Then $G$ will
   act freely on $X=S(V)$. These groups have been characterized by
   Wolf \cite{Wolf}:

   \begin{thm}
     The group $G$ will act freely on some $S(V)$ if and only if
     \begin{enumerate}
     \item every subgroup of order $pq$ ($p$, $q$ prime) in $G$ is
       cyclic, and
     \item $G$ does not contain $SL_2(\Fp)$ with $p>5$ as a subgroup.
     \end{enumerate}
   \end{thm}
 \end{exs}

\begin{NumberedSubSection}{General restrictions}

Consider a closed oriented manifold $M^n$ with a free $G$-action preserving
orientation. Then $M\to M/G$ has \emph{degree} $|G|$. Using the
spectral sequence for the fibration $M\to M/G\to BG$, one can prove

\begin{thm}
  (Browder \cite{Browder})
  If $G$ acts freely on $M^n$, preserving orientation,
  then $|G|$ divides the product
 \[ \prod_{s=1}^n\operatorname{exponent}
 H^{s+1}(G,H^{n-s}(M,\Z))\,.
\]
\end{thm}
\noindent In the statement above, note that the cohomology
groups which appear
are all finite abelian groups, and the \textsl{exponent}
is simply the smallest positive integer
that annihilates every element in the group.

\medskip

\begin{titled}{Consequences}
(1) If $M^n$ has a free action of $(\Z/p)^r$ which is trivial in
homology, then the total number of dimensions $0\le j< n$ such that
$H^j(M,\Z_{(p)})\ne 0$ must be at least~$r$. This follows from the
fact that $p\cdot \bar H^*((\Z/p)^r,\Z)=0$.
\par\smallskip\noindent(2)
If $(\Z/p)^r$ acts freely and homologically trivially on
$M=(S^n)^k$, then $r\le k$ (Carlsson \cite{Ca1}, \cite{Ca2}).
\end{titled}

With the help of Tate cohomology, these results can be extended to
finite, connected, free $G$--CW complexes. The most general
conjecture is given by

\begin{conj}
  (Carlsson \cite{Ca3}, \cite{Ca4}) If $(\Z/p)^r$ acts freely on $X$, then
  \[
       2^r \le \sum_{i=0}^{\dim X} \dim_{\Fp} H_i(X,\Fp)\,.
   \]
\end{conj}

\begin{titled}{Ancient open problem}
If $(\Z/p)^r$ acts freely on $X\simeq S^{n_1}\times\cdots S^{n_k}$,
then is $r\le k$?
\end{titled}

Let
\[
\begin{aligned}
       r(G)&=\max_{p\vert |G|} r_p(G)\\
       h(G) &=\min\left\{ k\,\vert\, G \text{ acts freely on a finite }X\simeq
              S^{n_1}\times\cdots\times S^{n_k}\right\}\,.
\end{aligned}
\]
then the following conjecture formulated by Benson and Carlson
(see \cite{BC}) is rather intriguing:

 \CenterBox{ \emph{Conjecture}:
  $r(G)=h(G)$}

\medskip
Swan's Theorem implies this result for $r(G)=1$. Note that
\emph{every} $G$ acts freely on \emph{some} product of spheres, so
$h(G)$ is well-defined.
\end{NumberedSubSection}

\begin{NumberedSubSection}{Case study : $S^n\times S^m$}

\begin{titled}{Problem}
Characterize those finite groups which act freely on a product of two
spheres.
\end{titled}

We will use representations and bundle theory to address this
problem, following the work of Adem, Davis and \"Unl\"u
\cite{ADU}.

\begin{defn}
  Let $G\subset U(n)$. We define the \emph{fixity} of this
  representation as the smallest integer $f$ such that $G$ acts freely
  on
$$U(n)/U(n-f-1).$$
\end{defn}

Note that $G$ has fixity zero if and only if $G$ acts freely on
$S^{2n-1}$, or equivalently if and only if $S^{2n-1}/G$ is a (complex)
linear space form.

There is a fibration sequence
\[
U(n-f)/U(n-f-1) \to U(n)/U(n-f-1) \to U(n)/U(n-f)
\]
in which the map from the total space to the base is a $G$-map.

\begin{thm}
  If $G\subset U(n)$ has fixity one, then $G$ acts freely and smoothly
  on $X=S^{2n-1}\times S^{4n-5}$.
\end{thm}

\begin{proof}
  Consider the bundle
\[
       U(n-1)/U(n-2)\to U(n)/U(n-2)\to U(n)/U(n-1)\,.
\]
It is the associated spherical bundle of a $G$-vector bundle $\xi$
such that $\xi\oplus\xi$ is trivial. Hence $S(\xi\oplus\xi)$
splits (non--equivariantly) as the indicated product of spheres, and
has a free $G$-action.
\end{proof}

\begin{cor}
  If $G\subset SU(3)$ is a finite subgroup,
then $G$ acts freely and smoothly on $S^5\times
  S^7$; the finite subgroups of $SU(3)$ include $A_5$, $SL_3(\Ftwo)$,
  and $3\cdot A_6$.
\end{cor}

We now focus on $p$-groups.

\begin{thm}
  If $p\ge3$ is a prime, then a finite $p$-group $P$ acts freely and
  smoothly on some $S^n\times S^m$  if and only if $P$ does not
  contain $\Z/p\times\Z/p\times\Z/p$ as a subgroup.
\end{thm}

\begin{titled}{Sketch of proof}
The ``only if'' statement has been known for 50 years
\cite{Heller}. For $p>3$,
one can show that either $P$ acts freely on some $S(V)\times S(W)$
for some representations $V$ and $W$, or else that there is a
representation $P\subset U(p)$ of fixity one, whence $P$ acts
freely on $S^{2p-1}\times S^{4p-5}$. This result involves using a
detailed description of rank two $p$-groups.

Partial results in the case $p=2$ are due to \"Unl\"u \cite{Unlu}.
There are 396 $2$-groups $P$ of order dividing $256$ and such that
$r_2(P)=2$. Of these only \emph{one} is not yet known to act
freely and smoothly on some $S^n\times S^m$.
\end{titled}

More generally, we have the

\begin{thm}
  (The Propagation Theorem) Let $G\subset U(n)$ be such that $G$ acts
  freely on $U(n)/U(k)$ for some $k>1$. Then, if
$$(|G|, (n-1)!)=1,$$
  $G$ will act freely and smoothly on
  \[
          M= S^{2n-1}\times S^{2n-3}\times\cdots\times S^{2k+1}\,.
\]
\end{thm}

\begin{cor}
  Let $P$ be a finite $p$-group with
  \begin{itemize}
  \item cyclic center,
  \item a maximal subgroup which is abelian, and
  \item rank $f+1 <p$.
  \end{itemize}
 Then $P$ acts freely and smoothly on
 \[
       M = S^{2p-1}\times\cdots\times S^{2(p-f)-1}\,,
 \]
 i.e., on a product of $f+1$ spheres.
\end{cor}

\begin{ex}
  Suppose that $P$ is an extra-special $p$-group of order $p^3$ and
  exponent~$p$.
  Then $P$ acts freely and smoothly on $S^{2p-1}\times S^{2p-3}.$
\end{ex}

\end{NumberedSubSection}

\begin{NumberedSubSection}{Homotopy actions}

The preceding results are explicit, geometric examples of a more
general homotopy--theoretic construction. The key ingredient is
the notion of a space with \emph{periodic cohomology}, which we
now define. The results here are based on the paper by Adem and
Smith \cite{AS}.

\begin{defn}
  A space $X$ is said to \emph{have periodic cohomology} if there
  exists a cohomology class $\alpha\in H^*(X,\Z)$, $|\alpha|>0$, and
  an integer $d\ge0$ such that for all coefficient systems $M$ the cup
  product map
  \[
     \cup\alpha: H^n(X,M)\to H^{n+|\alpha|}(X,M)
  \]
  is an isomorphism for $n\ge d$.
\end{defn}

\begin{exs}
  (1) If $G$ is a finite group, then $BG$ has periodic cohomology if
  and only if every abelian subgroup in $G$ is cyclic. Indeed, in this
  case periodicity for trivial coefficients implies periodicity for
  all coefficients.
\par\noindent(2) If $\Gamma$ is a discrete group of virtually finite
cohomological dimension, then $B\Gamma$ has periodic cohomology if and
only if every finite subgroup of $\Gamma$ has periodic cohomology.
\end{exs}

\begin{thm}
  (Adem--Smith) A connected CW--complex $X$ has periodic cohomology if
  and only if there exists an orientable spherical fibration
  \[
      S^N\to E\to X
  \]
  such that $E$ has the homotopy type of a finite dimensional complex.
\end{thm}

This has the following result as a corollary.

\begin{cor}
  If $Y$ is a simply--connected,
finite $G$--CW complex such that all of its isotropy
  subgroups have periodic cohomology, then there exists a finite
  \emph{free} $G$--CW complex
\[
    X\simeq S^N\times Y
\]
for some $N>\!\!> 0$.
\end{cor}

\begin{rems}
  Given $G\subset U(n)$ of fixity $f$, we have a geometric realization
  of our spherical fibration
 \[
     U(n-f)/U(n-f-1) \to U(n)/U(n-f-1)\to U(n)/U(n-f)\,.
  \]
\end{rems}

Using the corollary, we have reduced the problem of constructing free
$G$-actions on a finite complex $X\simeq S^n\times S^m$ to the problem
of \emph{constructing an action of $G$ on a sphere with rank one
  isotropy subgroups}.

\begin{exs}
  (1) Every rank two $p$-group $P$ acts on an $S(V)$ with rank one
  isotropy, hence $P$ acts freely on some finite $X\simeq S^n\times
  S^m$.
  \par\smallskip\noindent(2) If $S$ is a simple group of rank equal to two
  different from $PSL_3(\Fp)$, then $S$ acts freely on some finite
  $X\simeq S^n\times S^m$.
\end{exs}

Let $T_p= (\Z/p\times \Z/p)\rtimes SL_2(\Fp)$. Then one can show
that if $T_p$ acts on $X=S^n$, and $p$ is odd, it has a rank
\emph{two} isotropy subgroup. The equivariant cohomology
$H^*(ET_p\times_{T_p}X, \mathbb F_p)$ has Krull dimension equal to
two.

\begin{titled}{Problem}
Does $T_p$ act freely on some $S^n\times S^m$?
\end{titled}

Recently, M. Jackson \cite{Jackson} has announced

\begin{thm}
  A rank 2 group $G$ acts on $Y\simeq S^N$ with rank one isotropy if
  and only if $G$ does not contain $T_p$ as a subgroup for any odd
  prime~$p$. Consequently, all such groups act freely on a finite
  complex $X\simeq S^n\times S^m$.
\end{thm}

\noindent Hence we conclude that if $G$ is a rank~2 group not
containing any $T_p$ as a subgroup ($p$ odd), then $r(G)=h(G)=2$.

We will provide a more direct proof of this statement
for groups
of odd order by using some elementary group theory. Applying the
local methods in \cite{AS}, page 433, it suffices to prove
that for each prime $p$ such that $r_p(G)=2$, $G$ acts on a
sphere such that the isotropy subgroups have $p$--rank equal
to one.

\begin{thm}
If $G$ is a finite group of odd order and its rank is equal to
two, then for every prime $p$ such that $r_p(G)=2$, there exists a
$G$--representation $W_p$ such that the action of $G$ on $S(W_p)$
has isotropy subgroups having $p$--rank equal to one.
\end{thm}

\begin{proof}
First we need some group theory. It is known that every rank two
group of odd order is solvable and has a nilpotent commutator
subgroup (see \cite{Suzuki}, page 371). From this it follows
easily that for any prime $p$ dividing $|G|$, there exists a $p'$
normal subgroup $N$ such that $G'=G/N$ has a normal $p$--Sylow
subgroup.

Now let $G_p$ denote $Syl_p(G)=Syl_p(G')$, and assume that
$r_p(G)=2$. As in Theorem 3.8, we can find a representation $V_p$
for $G_p$ such that the $G_p$ action has rank one isotropy
(indeed, a central element of order $p$ acts freely on $V_p$ by
construction). Now we can induce this representation up to $G'$,
to obtain $W_p$; the associated
sphere $S(W_p)$ will have an action of $G'$ and hence of $G$ via
$G\to G'$ such that the isotropy subgroups have $p$-rank equal to
one.

\end{proof}

Hence we have

\begin{cor}
If $G$ is an odd order finite group of rank equal to two, then it
acts freely on a finite complex $X\simeq S^n\times S^m$.

\end{cor}

\end{NumberedSubSection}

\end{document}